%%%%%%%%%%%%%%%%%%%%%%%%%%%%%%%%%%%%%%%%%%%%%%%%%%%%%%%%%%%%%%%%%%%%%%%%%%%
%% Bryant, J.; Ferry, S.; Mio, W.; Weinberger, S.
%% 
%% Topology of homology manifolds
%% 
%% We construct examples of nonresolvable generalized $n$-manifolds, $n\geq
%%   6$, with arbitrary resolution obstruction, homotopy equivalent to any 
%%   simply connected, closed $n$-manifold. We further investigate the 
%%   structure of generalized manifolds and present a program for 
%%   understanding their topology.
%% 
%% publ:  Bull. Amer. Math. Soc. (N.S.) 28(1993) no. 2
%% pp:    324-328
%% type:  Research Announcement        markup: amstex    file size: 18K
%% contact:steve@math.binghamton.edu
%% 
%% copyright: American Math. Society copyright; see end of article
%% 
%% Include files necessary for this article: bull-ppt.tex
%% 
%%%%%%%%%%%%%%%%%%%%%%%%%%%%%%%%%%%%%%%%%%%%%%%%%%%%%%%%%%%%%%%%%%%%%%%%%%%
\input amstex
\documentstyle{amsppt}
\input bull-ppt
%\NoBlackBoxes
%\define\ni{\noindent}
\define\bZ{\Bbb Z}

\define\bR{\Bbb R}

%\define\qed{\vrule height 3pt width3pt depth2pt} 
%\newnumstyle\list1{\roman}

%\RefWarnings
%\baselineskip 12pt

%\def\mk{\medbreak}

%\hsize=20pc
\topmatter
\cvol{28}
\cvolyear{1993}
\cmonth{April}
\cyear{1993}
\cvolno{2}
\cpgs{324-328}
\title Topology of homology manifolds\endtitle 
\author J. Bryant, S. Ferry, W. Mio, and S. 
Weinberger\endauthor 
\address \RM{(J. Bryant and W. Mio)} Department of
Mathematics, Florida State University,
Tallahassee, Florida 32306\endaddress
\ml bryant\@math.fsu.edu
\mlpar {\it E-mail address\/}: mio\@math.fsu.edu \endml
\address \RM{(S. Ferry)} Department of Mathematics,
State University of New York at Binghamton, 
Binghamton, New York 13901\endaddress
\ml steve\@math.binghamton.edu\endml
\address \RM{(S. Weinberger)}
Department of Mathematics, University of Chicago, Chicago, 
Illinois
60637\endaddress
\thanks S. Ferry and S. Weinberger were partially 
supported by NSF 
grants. S. Weinberger was partially supported by a 
Presidential Young 
Investigator Fellowship\endthanks
%{\bf Abstract}.
\abstract We construct examples of nonresolvable 
generalized $n$-manifolds, $n\geq 6$, with arbitrary 
resolution 
obstruction, homotopy equivalent to any simply connected, 
closed 
$n$-manifold. We further investigate the structure of 
generalized 
manifolds and present a program for understanding their 
topology.
\endabstract
%\bigskip} 
\endtopmatter

\document
%\parindent=27pt
%{\narrower\smallskip
%\ni %\parindent=18pt
By a {\it generalized n-manifold} we will mean a 
finite-dimensional absolute
neighborhood retract $X$ such that $X$ is a homology 
$n$-manifold;
that is, for all $x\in X$, $H_i(X, X-\{x\}) =
H_i(\bR^n, \bR^n-\{0\})$. 
Generalized manifolds arise naturally as fixed-point sets 
of group actions
on manifolds,  as limits of sequences of manifolds, and as 
boundaries of negatively curved groups.
See [BM, Bo, B, GPW].
Such spaces have most of the homological properties of 
topological 
manifolds. In
particular, generalized manifolds satisfy Poincar\'e 
duality [Bo].
%\mk

Generalized manifolds also share certain geometric and 
analytic 
properties with
manifolds. Modern proofs of the topological invariance of 
rational 
Pontrjagin
classes show that Pontrjagin classes can be defined for 
generalized 
manifolds and
(even better!) that the symbol of the signature operator 
can be 
defined for
these spaces. See [CSW]. In light of this, the following 
question 
seems
natural:
%\mk
%\ni {\bf
\dfn{Question 1} Is every generalized manifold $X$ homotopy 
equivalent to a topological
manifold?
\enddfn
%\mk

By [FP], this is true for compact simply connected homology 
manifolds in all
higher dimensions. We shall see below that this is not 
true in the 
nonsimply
connected case. To continue in this vein, we can consider 
a strong 
version of
Question 1 that asserts that, for such an $X$\<, a 
manifold $M$ can be 
chosen
coherently for all of its open subsets.
%\mk
\dfn{Resolution conjecture \rm(see [C])} For every 
generalized 
$n$-manifold $X$
there is an $n$-manifold $M$ and a map $f:M\to X$ such 
that for 
each open $U$
in $X$, $f^{-1}(U)\to U$ is a (proper) homotopy 
equivalence.
\enddfn
%\mk

Quinn [Q2] showed that such a resolution is unique if it 
exists and 
reduced the
resolution conjecture for $X$ to the calculation of a 
locally defined 
invariant
$I(X) \in H^0(X; \bZ)$. For $X$ (and $Y$) connected, this 
invariant 
satisfies $I(X)\equiv 1$ (mod $8$), $I(X)=1$ iff $X$ has a 
resolution, 
and $I(X\times Y)=I(X)I(Y)$. He also showed that $I(X)$ is 
an
$s$-cobordism invariant
of generalized manifolds and that an affirmative solution 
to Question 
1 would imply
the resolution conjecture as well. We will henceforth 
assume that 
$X$ is
connected so that the obstruction $I(X)$ is an integer. %\mk

Resolutions are useful maps for studying the geometry of 
generalized manifolds.
For $n \ge 5$, Edwards [E] has characterized $n$-manifolds 
topologically as being
resolvable generalized manifolds that satisfy the 
following weak 
transversality
condition:
%\mk
%\ni{\bf 
\dfn{Disjoint disks property \rm(DDP)} $X$ has DDP if for 
any maps 
$f,g: D^2 \to
X$ and $\varepsilon>0$ there are maps $f',g':D^2 \to X$ 
with 
$d(f,f')<\varepsilon$,
$d(g,g')<\varepsilon$, and $f'(D^2) \cap
g'(D^2)=\emptyset$.
The resolution conjecture would then imply a 
characterization of 
topological
manifolds as generalized manifolds satisfying DDP. %\mk
\enddfn

Our first theorem says that the resolution conjecture is 
false.
%\mk
\proclaim{Theorem 1} Generalized
manifolds $X$ with arbitrary index
$I(X)\equiv 1\ (\roman{mod}\ 8)$ exist within the homotopy 
type of any 
simply connected
closed n-manifold, $n \ge 6$.\endproclaim %\mk
\proclaim{Corollary 1}{}There are generalized manifolds 
that are not 
homotopy equivalent to manifolds.\endproclaim %\mk
The remainder of this announcement is devoted to more 
positive 
statements. Some
are theorems and some are conjectures.
%\mk

The first statement requires a little preparation. Recall 
that an 
$s$-cobordism
is a manifold with two boundary components, each of which 
includes 
 a simple
homotopy equivalence. A basic result in
high-dimensional topology asserts that
$s$-cobordisms are products.
%\mk
\proclaim{Theorem 2} Let $\Cal S (M)$ denote the set of 
$s$-cobordism classes
of generalized manifolds mapping to $M$ by homotopy 
equivalences 
which are
homeomorphisms on the boundary. If $\dim\ M \ge 6$, then 
$\Cal S 
(M) = \Cal S (M
\times D^4)$.\endproclaim
%\mk	

This periodicity theorem was proven by Siebenmann [KS] for 
manifold structure sets,
but there was a $\bZ$ obstruction to its universal 
validity (see [N]). 
As suggested
by  Cappell, the theorem is valid in the class of 
generalized 
manifolds. The
following theorem will not
be precisely stated here but explains the functorial 
significance of 
these
spaces:
%\mk
\proclaim{Theorem 3} The algebraic surgery exact sequence 
{\rm(see [R])} is 
valid for
high-dimensional generalized manifolds up to 
$s$-cobordism.\endproclaim
%\mk

Thus, one knows that if $X$ is simply connected, $\Cal S 
(X)$ contains 
generalized
manifolds of every index and  there is a one-to-one 
correspondence 
between
$s$-cobordism classes of generalized manifolds homotopy 
equivalent 
to $X$ with
different indices. Simply connected surgery theory 
describes the 
structure of
the set of such generalized manifolds with $I(X)=1$. %\mk

On the other hand, for ``rigid" manifolds X like tori or 
nonpositively curved
manifolds, this theorem asserts that all generalized 
manifolds 
homotopy equivalent
to X are $s$-cobordant to the standard model. (This can be 
proven 
more directly
by using [FJ] and Quinn's resolution theorem [Q1, Q2]). %\mk

Our main conjectures are the following: %\mk
\dfn{Homogeneity conjecture} If $X$ is a connected 
generalized 
$n$-manifold
with DDP, $n \ge 5$, then given $p,q \in X$, there is a 
homeomorphism $h:X \to
X$ with $h(p)=q$. \enddfn
%\mk
\dfn{$S$-cobordism conjecture} $S$-cobordisms of DDP 
generalized manifolds are
products. \enddfn
%\mk

A corollary of these conjectures is that $n$-dimensional 
generalized manifolds of a given index that satisfy the 
DDP are 
``noncartesian manifolds" modeled on unique
local models. The following is very attractive and seems 
to be a step 
en route
to the $s$-cobordism conjecture.
%\mk
\dfn{Revised resolution conjecture} Every generalized 
manifold 
$X$ has a
resolution by a DDP generalized manifold. If $X$ satisfies 
DDP, then 
any such
resolution of $X$ is a uniform limit of
homeomorphisms.\enddfn
%\mk

With these ideas in place, we suggest a modification to 
another 
standard
conjecture:
%\mk
\dfn{Rigidity conjecture} If X is an aspherical Poincar\'e 
complex, then $X$ is
homotopy equivalent to a unique generalized manifold 
satisfying 
DDP.\enddfn
%\mk

If the polyhedron $X$ is nonpositively curved in the sense 
of [Gr], 
then one
can show, by using [FW], that the resolution obstruction 
is a well-defined
invariant of $\pi_1(X)$.
Farrell has informed the authors that 
this also
follows from the work of  Hu. The Borel rigidity 
conjecture (which for
nonpositively curved $\pi_1(X)$ also has been claimed by 
Hu) for 
manifolds with
boundary and the authors' realization theorem imply that a 
generalized 
manifold homotopy
equivalent to $X$ exists. The $s$-cobordism conjecture 
would then 
provide the
uniqueness. The authors' methods show that metric 
analogues of the usual 
rigidity
conjecture are false in situations where this version is 
correct (up to
$s$-cobordism).
%\mk

Construction of a counterexample among groups of 
nonpositive 
curvature could be
very useful in terms of producing natural examples of the 
anticipated local
models: these would arise as the spaces at infinity in 
natural 
compactifications
of the group (see [Gr]).
%\mk

Finally, we close with the following problem: 
\dfn{Dimension conjecture} The local index described above 
is 
realized in dimension four, but no three-dimensional 
examples of 
this sort exist.\enddfn
%\mk
%\ni {\smc 
\rem{Remarks on the proof} The spaces are built using 
controlled surgery as
developed in [Q1], [Q2], and [FP]. In surgery one tries to 
take 
manifold
approximations to a space and make them homotopy 
equivalent to 
the target. Our
construction is a ``controlled" or local version of this. 
The first step in 
the
construction is to construct a model of a Poincar\'e space 
$X$ by 
gluing together
two manifolds with boundary by using a homotopy equivalence 
between their
boundaries.
\endrem
%\mk

Successive approximations are refinements of this first 
step. 
We construct a
model for $X$ with better local Poincar\'e duality by 
gluing together 
two manifold
pieces with a homotopy equivalence which is controlled 
over the 
first stage.
The common boundary is constructed in such a way that its 
image in 
the first
stage is nearly dense. To produce the required controlled 
homotopy 
equivalence, 
controlled surgery theory comes into the picture. The new 
generalized
manifolds ``predicted'' by the surgery theory result from 
the 
difference between
controlled and uncontrolled surgery obstructions on the 
putative $X$.
%\mk

Each stage in the construction is a space with better 
Poincar\'e 
duality measured
over the previous stage of the construction. As the 
control on the 
gluing
homotopy equivalences improves, the resulting spaces have 
better 
local Poincar\'e
duality. This forces the limit of the approximations to be 
a homology 
manifold
because $X$ is a homology manifold if and only if $X$ 
satisfies 
Poincar\'e duality
locally measured over itself, that is, if and only if the 
constant sheaf 
is
Verdier dual to itself in the derived category of $X$. %\mk

Homotopy equivalences with good metric
properties are precisely the output of
the surgery theory of [FP]. To achieve the desired goal, 
the 
codimension-one
submanifold we glue along will become denser and denser in 
$X$. 
This is
inevitable, since the resolution obstruction can be 
measured on any 
open subset
of a homology manifold $X$. The nonresolvable generalized 
manifold 
is obtained as an
inverse limit of these approximations. In the limit, all 
of the 
approximate 
self-dualities become a genuine local self-duality. %\mk

The surgery theory for generalized manifolds follows from 
more controlled and
relative versions of the basic construction. 

%\bigskip
%\ni{\bf References.}
%\mk

\enddocument